\documentclass{article}
\usepackage{amssymb,amsmath}
\usepackage{graphics}

\def\qed{\hfill {\hbox{${\vcenter{\vbox{               
   \hrule height 0.4pt\hbox{\vrule width 0.4pt height 6pt
   \kern5pt\vrule width 0.4pt}\hrule height 0.4pt}}}$}}}

\def\tr{\triangleright}

\newtheorem{theorem}{Theorem}
\newtheorem{definition}{Definition}

\newtheorem{corollary}[theorem]{Corollary}
\newtheorem{example}{Example}

\newenvironment{proof}[1][Proof]{\smallskip\noindent{\bf #1.}\quad}%
{\qed\par\medskip}

\author{{\begin{tabular}{c} Natasha Harrell \\
\small{\texttt{natasha@math.ucr.edu}}\end{tabular}}
\and
{\begin{tabular}{c} Sam Nelson \\
\small{\texttt{knots@esotericka.org}}\end{tabular}}
\and
\small{\begin{tabular}{c}
Department of Mathematics, University of California, Riverside\\
900 University Avenue, Riverside, CA, 92521\end{tabular} }}

\date{}

\title{\Large \textbf{Non-classicality and quandle difference 
invariants}}

\begin{document}
\maketitle

\begin{abstract}
Non-classical virtual knots may have upper and lower quandles which are not 
isomorphic. We exploit this property to define quandle
difference invariants, which can detect non-classicality by comparing the 
numbers of homomorphisms into a finite quandle from a virtual knot's upper 
and lower quandles. The invariants for small-order finite quandles detect 
non-classicality in several interesting virtual knots. We compute the 
difference invariant with the six smallest connected quandles for all 
non-evenly intersticed Gauss codes with 4 crossings. For non-evenly 
intersticed Gauss codes with 4 crossings, the difference invariants detect
non-classicality in 86\% of codes which have non-trivial upper or lower 
counting invariant values.
\end{abstract}



\section{\large \textbf{Introduction}}

\textit{Virtual knots} are equivalence classes of Gauss codes under the
equivalence relation determined by Reidemeister moves. A Gauss code is
\textit{realizable} if it determines a planar knot diagram; virtual knots
which include realizable Gauss codes are \textit{classical} and those which 
do not are \textit{non-classical.} 
See \cite{K} and \cite{C1} for more. 

Detecting the non-classicality of a virtual knot is not always simple since 
each classical knot type includes some Gauss codes which cannot be realized 
without virtual crossings. Various methods have been proposed for detecting 
non-classicality of Gauss codes; see \cite{K2} and \cite{D}.

In this paper we define \textit{quandle difference invariants} which can
detect non-classicality of certain virtual knots. 
We show that quandle difference invariants detect non-classicality in some 
interesting examples of virtual knots, including one virtual knot whose 
non-classicality is not detected by the (upper) knot quandle or by the
one or two strand bracket polynomials, though its non-classicality is
detected with the three-strand bracket polynomial in \cite{D}. We describe an 
algorithm and provide an implementation in Maple for computing these 
invariants from a Gauss code and a finite quandle matrix.

In order to test the effectiveness of quandle difference invariants at 
detecting non-classicality, we generated a list of all Gauss codes with 
three and four crossings, eliminating codes which are obviously classical or 
trivial. From this list, we then computed the quandle difference invariant 
for each such code; the results are collected in table \ref{t1}.

The paper is organized as follows. Section 2 gives a brief review of virtual 
knots, section 3 introduces the quandle difference invariants and 2-component 
invariants, and section 4 contains computational results. The \texttt{Maple} 
code used to obtain these results is available at 
\texttt{www.esotericka.org/quandles}; bugfixes and improvements will be 
made as necessary.

\section{\large \textbf{Virtual knots and non-classicality}}

We begin with a definition from \cite{K}.

\begin{definition}\textup{
Let $K$ be an oriented knot diagram with crossings labeled with names and 
signs. A \textit{Gauss code} is a sequence of crossing labels with over/under
and sign information recorded in the order encountered while following the
orientation of the knot; it is defined up to cyclic permutation. A 
\textit{Gauss diagram} is obtained from a Gauss code by writing the code 
counterclockwise around a circle and joining the two instances of each crossing
with an arrow oriented ``in the direction of gravity'', that is, toward the 
undercrossing label, and noting the sign for each such arrow.}
\end{definition}

A Gauss code or diagram determines a knot diagram in a neighborhood of each 
crossing and specifies which strands are to be connected; thus a knot diagram 
can be reconstructed from its Gauss code. We can then consider doing knot 
theory combinatorially by defining a ``knot'' to be an equivalence class of 
Gauss codes under the equivalence relation determined by the Reidemeister 
moves. However, there is no guarantee that a given Gauss code determines a 
planar knot diagram; we may be forced to introduce additional crossings not 
mentioned in the code in order to draw the diagram in the plane. These 
\textit{virtual crossings} are drawn as circled self-intersections, and 
since they're not really there, we're permitted to redraw any section of an 
arc with only virtual crossings on its interior arbitrarily provided we put
only virtual crossings in the interior.

\begin{definition}\textup{
A \textit{virtual knot} is an equivalence class of Gauss codes under the
equivalence relation generated by the Reidemeister moves.}
\end{definition}

\begin{example} \label{e1} \textup{
The Gauss code \[UA+OB-UC+OD+OA+UB-UD+OC+\] corresponds to the pictured
virtual knot diagram and Gauss diagram.}
\[ \includegraphics{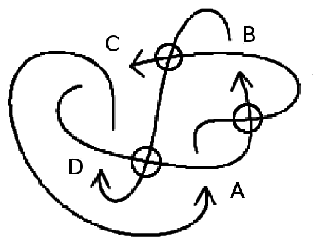} \quad \includegraphics{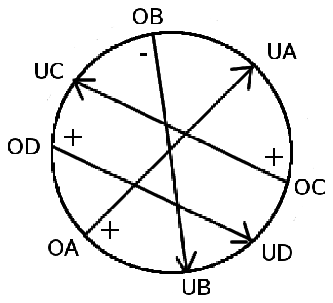} \]
\end{example}

Virtual crossings may be avoided by drawing our non-planar knot diagrams
on a surface with genus; we may then regard the Gauss code as describing a
knot in a thickened surface $\Sigma \times [0,1]$. The equivalence relation
defining virtual knots then translates into isotopy of the knot in 
$\Sigma\times [0,1]$ modulo stabilization and destabilization moves 
\cite{C1}, \cite{Ku}.

Classical knot theory forms a subset of virtual knot theory, since every 
classical knot may be considered a virtual knot without virtual crossings,
and, crucially, virtually isotopic classical knots are classically isotopic
\cite{PGV}, \cite{Ku}, \cite{N1}. Many invariants of classical knots extend 
to virtual knots by simply ignoring virtual crossings; these include the 
knot group, the knot quandle and biquandle, the various skein polynomials,
finite type invariants, etc.

Given a Gauss code, one asks whether the Gauss code corresponds to a classical
knot. If reconstruction yields a planar knot diagram, then the answer is 
clearly yes; however, every classical knot has representative Gauss codes 
which determine non-planar, i.e., virtual, diagrams, \cite{N1}. Thus, if 
a Gauss code is not obviously classical, the code may well represent either 
a classical or non-classical virtual knot. In the next section, we discuss
a method for detecting non-classicality of virtual knots using quandle 
difference invariants.

\section{\large \textbf{Quandle difference invariants}}

\begin{definition}
\textup{A \textit{quandle} is a set $Q$ with an operation 
$\tr:Q\times Q\to Q$ satisfying
\newcounter{q}
\begin{list}{(\roman{q})}{\usecounter{q}}
\item{for all $a\in Q$, $a\tr a = a$,}
\item{for every $a,b\in Q$ there is a unique $c\in Q$ such that $a=c\tr b$, 
and}
\item{for all $a,b,c \in Q$, $(a\tr b)\tr c=(a\tr c)\tr(b\tr c).$}
\end{list}
}\end{definition}

The term \textit{quandle} is due to Joyce; in \cite{J}, a quandle is 
associated to each knot diagram by a Wirtinger-type presentation. By 
construction, this knot quandle is an invariant of Reidemeister moves
and hence virtual isotopy. While it is difficult to directly compare 
quandles given by presentations, we can exploit the fact that knot quandles
are finitely presented to obtain a computable invariant of quandle 
isomorphism type by counting homomorphisms onto a finite quandle. 
This quandle counting invariant is studied in various papers and has been
jazzed up with quandle 2-cocyles from various quandle homology theories
\cite{D2}, \cite{C3}, \cite{C4}.

In \cite{PGV} it is observed that virtual knots have two knot groups,
the \textit{upper group} obtained by the usual Wirtinger presentation and
the \textit{lower group} obtained by first flipping the virtual knot over,
i.e., looking at the knot diagram from the other side of its supporting
surface, and then using the usual Wirtinger presentation. For classical
knots, flipping over is an ambient isotopy, so for classical knots, the
upper and lower groups are isomorphic. However, flipping over virtual 
knots involves doing forbidden moves which potentially change the virtual 
knot type (see \cite{N2}), so in general a virtual knot can have non-isomorphic
upper and lower groups. Exactly the same observation applies to quandles 
\cite{K}. 

\begin{theorem}
Let $K$ be a virtual knot diagram and let $\mathrm{Flip}(K)$ be the virtual
knot diagram obtained from $K$ by viewing $K$ from a base point on the other
side of the supporting surface. Then $\mathrm{Flip}(K)$ is an invariant of 
virtual isotopy.
\end{theorem}

\begin{proof}
Reidemeister and virtual Reidemeister moves on the flipped knot are simply 
Reidemeister and virtual Reidemeister moves on the original knot. 
\end{proof}

\begin{corollary} \label{th1}
Let $\Phi$ be any virtual knot invariant. Then 
\[\Phi_2(K)=(\Phi(K), \Phi(\mathrm{Flip}(K)))\] 
is an invariant of virtual isotopy.
\end{corollary}

The two-component version of a virtual knot invariant does not reveal any 
additional information about 
$K$ when $K$ is classical, but they can be very helpful in distinguishing
non-classical virtual knots. Obvious examples of 2-component invariants 
include the 2-component couting invariants from the knot group, quandle and 
biquandle; 2-component CJKLS 2-cocyle invariants; the 2-component Jones 
polynomial, etc.
 
Corollary \ref{th1} suggests the following

\begin{definition}\textup{
Let $T$ be a finite quandle and $K$ a virtual knot. Let $U$ denote the 
upper quandle of $K$ and $L$ the lower quandle of $K.$ The 
\textit{two-component counting invariant} is the ordered pair
\[Q2_T(K)=\left(|\mathrm{Hom}(U,T)|,|\mathrm{Hom}(L,T)|\right).\]
The 
\textit{quandle difference invariant} $QD_T(K)$ of $K$ associated to $T$
is the difference between the number of the number of homomorphisms
from $U$ to $T$ and the number of of homomorphisms from $L$ to $T$. That is,
\[QD_T(K)=|\mathrm{Hom}(U,T)| - |\mathrm{Hom}(L,T)|.\]
}
\end{definition}

The fact that classical knots have isomorphic upper and lower quandles implies

\begin{theorem}
If a virtual knot $K$ is classical, then $QD_T(K)=0$ for every finite quandle 
$T$.
\end{theorem}

\begin{corollary}
If $QD_T(K)\ne 0$ for any finite quandle $T$, then $K$ is not classical. In 
particular, if $QD_T(K)\ne 0$ for any finite quandle $T$, then $K$ is not 
the unknot.
\end{corollary}

Moreover, since the individual components of $QD_T$ are invariants of virtual 
isotopy, $QD_T$ is itself an invariant of virtual knots; $QD_T(K)\ne QD_T(K')$
implies $K$ is not virtually isotopic to $K'$. In the present work, however, 
we focus on using the quandle difference invariant to detect non-classicality 
of virtual knots. We know already that the quandle difference invariant does 
not always detect non-classicality.

\begin{example} \textup{
\[
\includegraphics{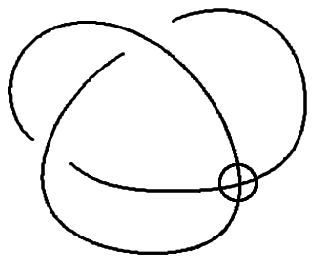} \qquad \includegraphics{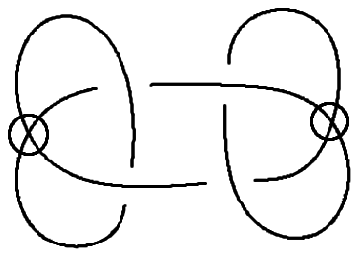}
\]
The knot on the left is Kauffman's virtual trefoil. It is known to be 
non-trivial and non-classical \cite{K}. However, its upper and lower quandles
are both isomorphic to the free quandle on one generator, that is, the knot
quandle of the unknot. The knot on the left is a variant of the Kishino knot; 
it has rotational symmetry about the horizontal axis, so its upper and lower 
quandles have identical presentations.
}\end{example}

\begin{definition}\textup{
We say a virtual knot is \textit{pseudoclassical} if it is virtually isotopic 
to its flip. In particular, a pseudoclassical virtual knot has isomorphic
upper and lower quandles.
}\end{definition}

Clearly, the quandle difference invariant cannot detect the non-classicality 
of pseudoclassical virtual knots for any finite quandle $T$. However, the
quandle difference invariants are effective at detecting
non-classicality in many virtual knots. Using the smallest connected quandle,
for example, we detect the non-classicality (and hence non-triviality) of the
non-alternating Kishino knot in figure 1. 
It only took two tries to find a quandle difference invariant which detects 
the non-classicality of the virtual knot $K_D$ in 1; this knot is 
not distinguished from the unknot by either the Jones polynomial or 2-strand 
bracket polynomial, though the 3-strand bracket polynomial does the trick 
\cite{D}.

\begin{figure}[!ht]\label{fig:kish}
\[K=\raisebox{-0.5in}{\includegraphics{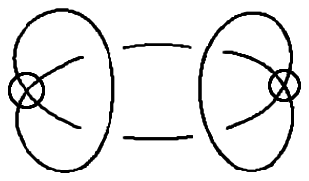}} 
\quad T=\left[\begin{array}{lll}
1 & 3 & 2 \\
3 & 2 & 1  \\
2 & 1 & 3 
\end{array}
\right]\]
\[QD_{T}(K)=-6\]
\[K_D=\raisebox{-0.5in}{\includegraphics{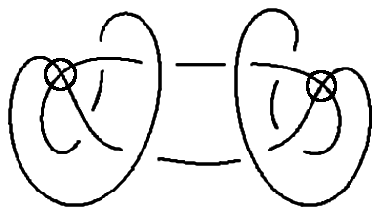}} 
\quad T_2=\left[\begin{array}{llll}
1 & 3 & 4 & 2 \\
4 & 2 & 1 & 3 \\
2 & 4 & 3 & 1 \\
3 & 1 & 2 & 4
\end{array}
\right]\]
\[QD_{T_2}(K_D)=-12\]
\caption{$QD_T$ detects non-classicality of virtual knots $K$ and $K_D$.}
\end{figure}

\begin{figure}[!hb]\label{gsigns}
\begin{center} 
\begin{tabular}{cc}
\raisebox{-0.6in}{\includegraphics{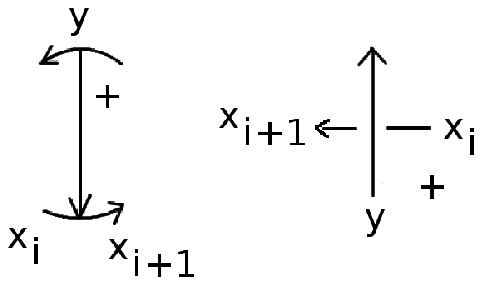}} &
$\Rightarrow x_i\tr y=x_{i+1}$ \\
\raisebox{-0.6in}{\includegraphics{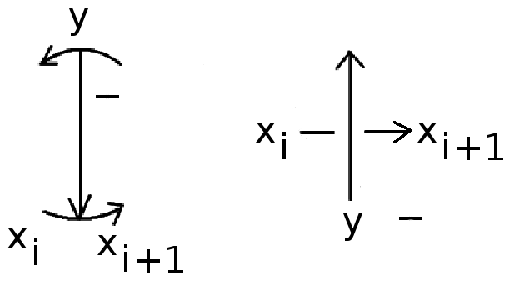}} &
$\Rightarrow x_{i+1}\tr y=x_i$
\end{tabular}
\end{center}
\caption{Quandle relations from a Gauss diagram.}
\end{figure}

To compute the upper and lower quandle presentations from a Gauss diagram, we 
divide the outer circle at each arrowhead to obtain a list of arcs in the 
diagram, and the quandle relations at each crossing are determined according 
to the crossing sign. As figure 2 shows, the rule is that if the arrowhead 
dividing $x_i$ from $x_{i+1}$ has tail on arc $y$, the relation determined is 
$x_i\tr y = x_{i+1}$ if the crossing sign  is  $+$  and $x_{i+1}\tr y = x_i$  
if the crossing sign  is $-$.

Dividing the same Gauss diagram at the arrowtails and applying the same 
procedure yields a presentation of the lower quandle $L$. Alternatively,
we can obtain a presentation of the lower quandle by reversing the direction
of the arrows, then getting the upper quandle presentation of the resulting
diagram.

\begin{example}\textup{ \label{e2}
Let $K$ be the virtual knot from example \ref{e1}. 
Then K has upper quandle presentation
\[U=\langle 1,2,3,4 \ | \ 1\tr 3 = 2, 2\tr 1 = 3, 4\tr 2 = 3, 4\tr 3 = 1 
\rangle\] and lower quandle presentation
\[
L=\langle 1,2,3,4 \ | \ 1\tr 3 = 2, 2\tr 4 = 3, 3\tr 1= 4, 1\tr 3 = 4\rangle.
\]
\[\includegraphics{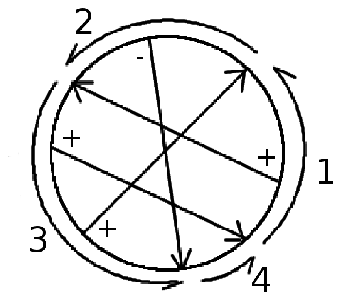} \quad \includegraphics{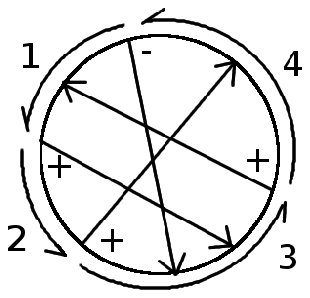}\]
}\end{example}

From these presentations, we can compute the quandle difference invariant
$QD_T(K)$ for any choice of finite target quandle $T$ by systematically 
considering all assignments of element of $T$ to each generator in the
quandle presentation and checking whether the assignment satisfies the
required relations in the target quandle. In the next section, we describe
our computational procedure for computing the quandle difference invariant
for any choice of Gauss code and finite target quandle and give some 
computational results.

\section{\large \textbf{Computational results}}

\texttt{Maple} code implementing the algorithms below is available in
the file \texttt{quandledifference.txt} at the second listed author's website
\texttt{www.esotericka.org/quandles}.

We compute the quandle difference invariant $QD_T(K)$ using the quandle 
matrix notation from \cite{HN}. Specifically, let $T=\{x_1\dots, x_n\}$ be 
a finite quandle. Then the matrix of $T$, $M_T$ is the matrix abstracted from 
the operation table of $T$ by dropping the $x$s,
keeping only the subscripts. A list of all finite quandles with up to 5 
elements appears in \cite{HN}; lists of all finite quandles with 6, 7 and 8 
elements along with software for computing such lists are available 
at \texttt{www.esotericka.org/quandles}.\footnote{An
independently derived list of quandles with up to 6 elements, obtained by 
Yamada, is included in \cite{C2}.} The file \texttt{quandles-maple.txt}
contains \texttt{Maple} code for counting homomorphisms from a knot quandle 
to finite quandle specified by a matrix.

We represent a finitely presented quandle with a \textit{quandle 
presentation matrix} as described in \cite{HMN}. That is, we put 
\[M_{ij}=\left\{\begin{array}{ll} k, & x_i\tr x_j = x_k, \\
0, & \mathrm{else}.
\end{array}\right.\]

If a Gauss code has a negative undercrossing label $x_i$ such that the next 
undercrossing label $x_{i+1}$ has sign $+$ and both overcrossings 
lie on the same arc $x_j$, then the quandle relations determined are 
$x_i\tr x_j=x_{i-1}$ and $x_i\tr x_j=x_{i+1}$, which determine conflicting 
entries for the $i,j$ position in the presentation matrix. In this situation 
we can do a type I Reidemeister move to obtain a presentation with no
conflicting entries; the program \texttt{gfix} finds and corrects this
problem.

\begin{example}\textup{
The virtual knot diagram in example \ref{e1} has upper and lower quandle 
presentation matrices
\[ U=\left[
\begin{array}{rrrr}
0 & 0 & 2 & 0 \\
3 & 0 & 0 & 0 \\
0 & 0 & 0 & 0 \\
0 & 3 & 1 & 0
\end{array}
\right] \quad \mathrm{and} \quad
L =\left[
\begin{array}{rrrrr}
0 & 0 & 2 & 0 & 0 \\
0 & 0 & 0 & 0 & 3 \\
0 & 0 & 0 & 0 & 0 \\
0 & 3 & 0 & 5 & 0 \\
0 & 1 & 0 & 0 & 0
\end{array}\right].
\]
}\end{example}

We represent Gauss codes in \texttt{Maple} as vectors with (appropriately?) 
Gaussian integer entries, using the convention
\[ OX+ \leftrightarrow X, \quad UX+ \leftrightarrow -X, \] 
\[
OX- \leftrightarrow X+I, \quad UX- \leftrightarrow -X-I, \quad 
X\in\mathbb{Z}. \]
To allow for multi-component links, each component is ended with a ``0''.

\begin{example}\textup{
The multicomponent Gauss code 
\[ U1+O2-O3-U2-, \quad O1+U3-\ \]
is represented in our \texttt{Maple} notation as the vector
\[
[-1,2+I,3+I,-2-I,0,1,-3-I,0].
\]
}\end{example}

Our quandle-difference computing algorithm then takes a Gauss code and
first runs it through \texttt{gfix}, which detects the situation described 
above and does Reidemeister I moves as necessary until the code has an 
upper quandle presentation matrix with no conflicting entries. The program 
\texttt{gauss2pres} then reads off the upper quandle presentation matrix from 
the Gauss code. Our main program, \texttt{qdiff}, uses \texttt{homcount} from 
\texttt{quandles-maple.txt} to compute the number of homomorphisms from 
the upper quandle of the matrix returned by \texttt{gfix} into the specified
target quandle. Then, the program simply multiplies each entry in the code by 
$-1$ to switch the direction of the arrows and cycles through to put the 
$-1$ or $-1-i$ entry in the first position, then repeats the procedure to get 
the lower quandle presentation matrix and its counting invariant. Finally,
the program reports the difference between these two numbers. We also include
\texttt{q2chom} which computes the 2-component quandle counting invariant.

In order to get some feeling for the effectiveness of quandle difference 
invariants at detecting non-classicality of virtual knots, we generated 
lists of all non-evenly intersticed single-component Gauss codes with three 
and four crossings,
removing codes which reduce by an obvious type I or type II move, using
\texttt{rglist}. The non-evenly intersticed condition guarantees 
non-planarity of the corresponding virtual knot, which may be classical or
non-classical (including pseudoclassical). There are 172 such 3-crossing Gauss 
codes; however, $QD_T$ does not detect non-classicality in any of these 
3-crossing
codes for the six smallest connected quandles. There are 17040 such 4-crossing 
codes. We then computed the number of codes in which non-classicality was 
detected using some connected quandles of order up to six. The results are 
collected in table 1. It is not clear what percentage of the Gauss codes with 
$QD_T=0$ are classical (i.e., diagrams of the unknot, trefoil, or figure 
eight), what percentage are pseudoclassical, and what percentage may have 
non-classicality detected by larger finite quandles, nor 
is it clear how many distinct virtual knots are represented. However, 
eliminating the codes which have trivial values for the quandle counting 
invariants for both upper and lower quandles for all six listed quandles, a 
step which eliminates the classical and pseduoclassical codes from the lists,
gives a list of 16 3-crossing codes and 4140 4-crossing codes. Of these 
4-crossing codes, 3570 or 86\% have non-classicality detected by at least one 
of the six listed quandle difference invariants.

\begin{table}[!ht] 
\begin{center} {\small
\begin{tabular}{|c|c|} \hline
Target quandle  &  
Number detected/ 17040 \\ 
$T$ & (4-crossing codes) \\ \hline
  &  \\
$\left[\begin{array}{lll}
1 & 3 & 2 \\
3 & 2 & 1 \\
2 & 1 & 3
\end{array}
\right]$  & 3060 \\  
 &   \\
$\left[\begin{array}{llll}
1 & 3 & 4 & 2 \\
4 & 2 & 1 & 3 \\
2 & 4 & 3 & 1 \\
3 & 1 & 2 & 4
\end{array}
\right]$   & 1350 \\ 
 &   \\
$\left[\begin{array}{lllll}
1 & 3 & 4 & 5 & 2 \\
3 & 2 & 5 & 1 & 4 \\
4 & 5 & 3 & 2 & 1 \\
5 & 1 & 2 & 4 & 3 \\
2 & 4 & 1 & 3 & 5
\end{array}
\right]$  &  492 \\ 
 &   \\
$\left[\begin{array}{lllll}
1 & 4 & 5 & 3 & 2 \\
3 & 2 & 4 & 5 & 1 \\
2 & 5 & 3 & 1 & 4 \\
5 & 1 & 2 & 4 & 3 \\
4 & 3 & 1 & 2 & 5
\end{array}
\right]$  &   72 \\ 
 &  \\
$\left[\begin{array}{lllll}
1 & 4 & 5 & 2 & 3 \\
3 & 2 & 1 & 5 & 4 \\
4 & 5 & 3 & 1 & 2 \\
5 & 3 & 2 & 4 & 1 \\
2 & 1 & 4 & 3 & 5
\end{array}
\right]$   &  426 \\ 
 &   \\
$\left[\begin{array}{llllll}
 1 & 4 & 5 & 2 & 3 & 1 \\
 4 & 2 & 6 & 1 & 2 & 3 \\
 5 & 6 & 3 & 3 & 1 & 2 \\
 2 & 1 & 4 & 4 & 6 & 5 \\
 3 & 5 & 1 & 6 & 5 & 4 \\
 6 & 3 & 2 & 5 & 4 & 6
\end{array}
\right]$  &  3060 \\
 &   \\ 
\hline
\end{tabular}}
\end{center}
\caption{Number of non-evenly intersticed single-component Gauss codes in 
which \texttt{qdiff} detects non-classicality.}
\label{t1}
\end{table}

\end{document}